\documentclass[reqno]{amsart}
\usepackage{amsmath}%
\usepackage{amsfonts}%
\usepackage{amssymb}%
\usepackage{graphicx}
\usepackage{setspace}
\usepackage{microtype}

\theoremstyle{plain}

\numberwithin{equation}{section}
\begin{document}
\onehalfspacing
\raggedbottom
\title[B-Splines]{Constructing Explicit B-Spline}
\author{R.O.~Linger}
\author{H.R.N.~van~Erp}
\author{P.H.A.J.M.~van~Gelder}

\begin{abstract}
We introduce here a direct method to construct multivariate explicit B-spline bases. B-splines are piecewise polynomials, which are defined on adjacent tetrahedra and which are $C^{r}$ continuous throughout. The $C^{r}$ continuity is enforced by making sure that all directional derivatives of order $r$, and lower, on the boundaries of adjacent tetrahedra give the same values for both tetrahedra. The method presented here is explicit, in that we will provide an algorithm with which one can analytically construct the B-spline base that enforces $C^{r}$ continuity for a given geometry. 
\end{abstract}
\maketitle

\section{Some Preliminaries}
In order to give a direct method for constructing B-splines, we must first give some preliminaries of these B-splines. In what follows, we will introduce the concepts of barycentric coordinates, barycentric polynomials, and directional derivatives on barycentric polynomials. This is done for two-dimensional functions.

Say, we have three points $\boldsymbol{v}_{1}$, $\boldsymbol{v}_{2}$, $\boldsymbol{v}_{3}$, in a two-dimensional Cartesian coordinate system with coordinates$\left(x,y\right)$. Then the difference vectors $\boldsymbol{v}_{1}-\boldsymbol{v}_{3}$ and $\boldsymbol{v}_{2}-\boldsymbol{v}_{3}$ are the axes of a coordinate system with origin $\boldsymbol{v}_{3}$ and coordinates, say, $b_{1}$ and $b_{2}$. We then have that
\begin{align}
	\label{eq.1}
	x = v_{13} + b_{1}\left(v_{11}-v_{13}\right) + b_{2}\left(v_{12}-v_{13}\right), \nonumber \\
	y = v_{23} + b_{1}\left(v_{21}-v_{23}\right) + b_{2}\left(v_{22}-v_{23}\right), 
\end{align}   
or, equivalently,
\begin{align}
	\label{eq.2}
	x = b_{1} v_{11}  + b_{2} v_{12} + \left(1 - b_{1}- b_{2}\right) v_{13},\nonumber \\
	y = b_{1} v_{21}  + b_{2} v_{22} + \left(1 - b_{1}- b_{2}\right) v_{23}.
\end{align}   

If a given coordinate $\left(x,y\right)$ lies in the triangle, say, $T$, with origin $\boldsymbol{v}_{3}$ and spanned by the vectors $\boldsymbol{v}_{1}-\boldsymbol{v}_{3}$ and $\boldsymbol{v}_{2}-\boldsymbol{v}_{3}$, then
\begin{equation}
	\label{eq.3}
	0 \leq b_{1} + b_{2} \leq 1.
\end{equation}
We define
\begin{equation}
	\label{eq.4}
	b_{3} \equiv 1 - b_{1} - b_{2}.
\end{equation}
By way of \eqref{eq.3} and \eqref{eq.4}, we have that the coordinates $\left(b_{1}, b_{2}, b_{3}\right)$, for a given point $\left(x,y\right)$ in triangle $T$, are all greater than zero. Furthermore, these coordinates sum to
\begin{equation}
	\label{eq.5}
	b_{1} + b_{2} + b_{3} = 1.
\end{equation}  

The coordinates $\left(b_{1}, b_{2}, b_{3}\right)$ are called barycentric coordinates. Using these coordinates, we may rewrite \eqref{eq.2} as
\begin{align}
	\label{eq.6}
	x = b_{1} v_{11}  + b_{2} v_{12} + b_{3} v_{13},\nonumber \\
	y = b_{1} v_{21}  + b_{2} v_{22} + b_{3} v_{23},
\end{align}   
or, equivalently,
\begin{equation}
	\label{eq.7}
	\left(\begin{array}{c} x\\y\end{array}\right) = \left[\begin{array}{ccc} \boldsymbol{v}_{1} & \boldsymbol{v}_{2} & \boldsymbol{v}_{3} \end{array}\right] \left(\begin{array}{c} b_{1}\\b_{2}\\b_{3} \end{array}\right).
\end{equation}   

Say, we have two barycentric coordinate systems $T_{1}$ and $T_{2}$:
\begin{figure}[!h]
	\centering
		\includegraphics[width=1.0\textwidth]{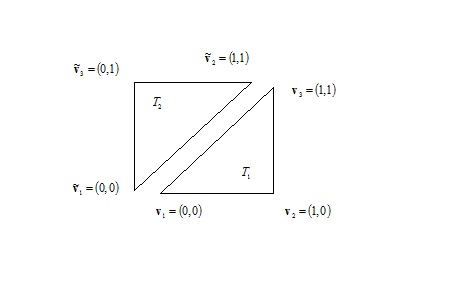}
	\caption{Simple Triangulation}
	\label{fig:triangulation}
\end{figure}
\\
\\ 
\noindent Then a given Cartesian point $\left(x,y\right)$ in either $T_{1}$ or $T_{2}$ may be written as, respectively,
\begin{equation}
	\label{eq.8}
	\left(\begin{array}{c} x\\y\end{array}\right) = \left[\begin{array}{ccc} \boldsymbol{v}_{1} & \boldsymbol{v}_{2} & \boldsymbol{v}_{3} \end{array}\right] \left(\begin{array}{c} b_{1}\\b_{2}\\b_{3} \end{array}\right), \qquad \left(x,y\right)\in T_{1},
\end{equation}   
and
\begin{equation}
	\label{eq.9}
	\left(\begin{array}{c} x\\y\end{array}\right) = \left[\begin{array}{ccc} \tilde{\boldsymbol{v}}_{1} & \tilde{\boldsymbol{v}}_{2} & \tilde{\boldsymbol{v}}_{3} \end{array}\right] \left(\begin{array}{c} \tilde{b}_{1}\\\tilde{b}_{2}\\\tilde{b}_{3} \end{array}\right), \qquad \left(x,y\right)\in T_{2},
\end{equation}   
We call the collection of triangles $T_{1}$ and $T_{2}$ a triangulation of the Cartesian plane. 

Note that for a given point $\left(x,y\right)$ in, say, $T_{1}$ the corresponding barycentric coordinates may be found by way of
\begin{equation}
	\label{eq.10}
	\left(\begin{array}{c} b_{1}\\b_{2}\end{array}\right) = \left[\begin{array}{cc} v_{11}-v_{13} & v_{12}-v_{13} \\ 
	v_{21}-v_{23} & v_{22}-v_{23}\end{array}\right]^{-1} \left(\begin{array}{c} x - v_{13}\\ y - v_{23} \end{array}\right),
\end{equation}   
and, by construction, \eqref{eq.4}:
\[	
	b_{3} = 1 - b_{1} - b_{2}.
\]

If we have some $d$th-order polynomial $p$ defined on the whole Cartesian plane $\left(x,y\right)$, say:
\begin{equation}
	\label{eq.11}
	p\!\left(x,y\right) = \sum_{0\leq i+j\leq d} \alpha_{ij}  x^{i}  y^{j}
\end{equation}  
Then we may make a Jacobian transformation from $\left(x,y\right)$ to $\left(b_{1}, b_{2}, b_{3}\right)$, by way of \eqref{eq.10} and \eqref{eq.4}, in order to obtain a polynomial in barycentric coordinates, equivalent to \eqref{eq.11}:
\begin{equation}
	\label{eq.12}
	p\!\left(b_{1},b_{2},b_{3}\right) = \sum_{0\leq i+j\leq d} \gamma_{ij} \: b_{1}^{i} \: b_{2}^{j} \: b_{3}^{d-i-j}.
\end{equation} 	
Likewise, we may make a Jacobian transformation from $\left(x,y\right)$ to $\left(\tilde{b}_{1}, \tilde{b}_{2}, \tilde{b}_{3}\right)$, in order to obtain the equivalent polynomial:
\begin{equation}
	\label{eq.13}
	p\!\left(\tilde{b}_{1},\tilde{b}_{2},\tilde{b}_{3}\right) = \sum_{0\leq i+j\leq d} \tilde{\gamma}_{ij} \: \tilde{b}_{1}^{i} \: \tilde{b}_{2}^{j} \: \tilde{b}_{3}^{d-i-j}.
\end{equation}   
Note that in \cite{Awanou03}, a more complicated proof of the equivalence of, say, \eqref{eq.11} and \eqref{eq.12} is given. But such a proof is unnecessary, in that a simple change of variable argument will suffice.

In what follows, we will have to make use of the fact that $r$th order continuity of two different piecewise polynomials defined, respectively, on two connected triangles $T_{1}$ and $T_{2}$, implies that all $r$th-order directional derivatives of these polynomials be equal on their shared boundary. So, we introduce here the concept of directional derivatives of polynomials in barycentric coordinates. 

Let $\boldsymbol{u}$ be the directional vector
\begin{equation}
	\label{eq.14}
	\boldsymbol{u} = a_{1} \boldsymbol{v}_{1} + a_{2} \boldsymbol{v}_{2} + a_{3} \boldsymbol{v}_{3}. 
\end{equation}   
Then the barycentric coordinates of $\boldsymbol{u}$ in $T_{1}$ are
\begin{equation}
	\label{eq.15}
	\boldsymbol{a} = \left(a_{1},a_{2},a_{3}\right)^{T}. 
\end{equation}   
Let $D_{\boldsymbol{u}}^{\left(r\right)}$ be the first-order directional derivative operator in the Cartesian plane $\left(x,y\right)$. Then we may define the zeroth-order directional derivative, in both the Cartesian and barycentric coordinate systems, as:
\begin{equation}
	\label{eq.16}
	D_{\boldsymbol{u}}^{\left(0\right)} p\!\left(x,y\right) = p\!\left(x,y\right) = p\!\left(b_{1},b_{2},b_{3}\right) = D_{\boldsymbol{a}}^{\left(0\right)} p\!\left(b_{1},b_{2},b_{3}\right).
\end{equation}   
The first-order directional derivative is defined as:
\begin{align}
	\label{eq.17}
	D_{\boldsymbol{u}}^{\left(1\right)} p\!\left(x,y\right) &= \boldsymbol{u}^{T} \nabla p\!\left(x,y\right)\nonumber \\
	\nonumber\\
	&= \boldsymbol{a}^{T} \nabla p\!\left(b_{1},b_{2},b_{3}\right) \\
	\nonumber\\
	&= D_{\boldsymbol{a}}^{\left(1\right)} p\!\left(b_{1},b_{2},b_{3}\right). \nonumber
\end{align}
where $\nabla$ is defined, respectively, as
\[
	\nabla = \left(\frac{\partial}{\partial x}, \frac{\partial}{\partial y}\right), \qquad \text{and}\qquad \nabla = \left(\frac{\partial}{\partial b_{1}}, \frac{\partial}{\partial b_{2}}, \frac{\partial}{\partial b_{3}}\right),
\]
depending on the coordinate system on which the polynomial $p$ is defined. By way of \eqref{eq.16} and \eqref{eq.17}, we may define the $r$th-order directional derivative recursively as:
\begin{align}
	\label{eq.18}
	D_{\boldsymbol{u}}^{\left(r\right)} p\!\left(x,y\right) &= D_{\boldsymbol{u}}^{\left(1\right)} D_{\boldsymbol{u}}^{\left(r-1\right)} p\!\left(x,y\right) \nonumber \\
	\nonumber\\
	&= D_{\boldsymbol{a}}^{\left(1\right)} D_{\boldsymbol{a}}^{\left(r-1\right)} p\!\left(b_{1},b_{2},b_{3}\right) \\
	\nonumber\\
	&= D_{\boldsymbol{a}}^{\left(r\right)} p\!\left(b_{1}, b_{2}, b_{3}\right). \nonumber
\end{align}   

For example, the zeroth-, first-, and second-order directional derivatives for arbitrary $\boldsymbol{a}$ and a second-order piecewise polynomial may be written down as, \eqref{eq.12} and \eqref{eq.16}:
\begin{align}
	\label{eq.19}
	 D_{\boldsymbol{a}}^{\left(0\right)} p\!\left(b_{1}, b_{2}, b_{3}\right) &= \sum_{0\leq i+j\leq 2} \gamma_{ij} \: b_{1}^{i} \: b_{2}^{j} \: b_{3}^{d-i-j} \nonumber \\
	 \\
	 &= \gamma_{00} \:  b_{3}^{2} + \gamma_{10} \:  b_{1}   b_{3} + \gamma_{01} \:  b_{2}   b_{3} + \gamma_{11} \:  b_{1}   b_{2} + \gamma_{20} \:  b_{1}^{2} + \gamma_{02} \:  b_{2}^{2}\nonumber 
\end{align}
and
\begin{align}
	\label{eq.20}
	 D_{\boldsymbol{a}}^{\left(1\right)} p\!\left(b_{1}, b_{2}, b_{3}\right) &= \boldsymbol{a}^{T} \nabla \sum_{0\leq i+j\leq 2} \gamma_{ij} \: b_{1}^{i} \: b_{2}^{j} \: b_{3}^{d-i-j} \nonumber \\
	 \nonumber \\
	 &= \boldsymbol{a}^{T} \left(
\begin{array}{c}
		\gamma_{10} \:  b_{3} + \gamma_{11} \:  b_{2}  +  2 \gamma_{20} \:  b_{1} \\
		\gamma_{01} \:  b_{3} + \gamma_{11} \:  b_{1}  +  2 \gamma_{02} \:  b_{2} \\
		2 \gamma_{00} \:  b_{3} + \gamma_{10} \:  b_{1} + \gamma_{01} \:  b_{2} 
\end{array}
\right) \\
	  \nonumber \\
	  &= \gamma_{00} \: 2 a_{3} b_{3} + \gamma_{10} \left(a_{1} b_{3}+a_{3} b_{1}\right) + \gamma_{01} \left(a_{2} b_{3}+a_{3} b_{2}\right) \nonumber \\ 
	  & \left.\qquad\right. + \gamma_{11} \left(a_{1} b_{2}+a_{2} b_{1}\right) +  \gamma_{20} \: 2 a_{1} b_{1} +  \gamma_{02} \: 2 a_{2} b_{2}.\nonumber 
\end{align}
and
\begin{align}
	\label{eq.21}
	 D_{\boldsymbol{a}}^{\left(2\right)} p\!\left(b_{1}, b_{2}, b_{3}\right) &= \boldsymbol{a}^{T} \nabla\left[ D_{\boldsymbol{a}}^{\left(1\right)} p\!\left(b_{1}, b_{2}, b_{3}\right) \right] \nonumber \\
	 \nonumber \\
	 &= \boldsymbol{a}^{T} \left(
\begin{array}{c}
	\gamma_{10} \: a_{3} +  \gamma_{11} \: a_{2} +  \gamma_{20} \: 2 a_{1} \\
	\gamma_{01} \: a_{3} +  \gamma_{11} \: a_{1} +  \gamma_{02} \: 2 a_{2} \\
	\gamma_{00} \: 2 a_{3} + \gamma_{10} \: a_{1} + \gamma_{01} \: a_{2} 
\end{array}
\right) \\
	  \nonumber \\
	  &=  \gamma_{00} \: 2 a_{3}^{2} + \gamma_{10} \: 2 a_{1} a_{3} + \gamma_{01} \: 2 a_{2} a_{3} \nonumber \\ 
	  & \left.\qquad\right.  + \gamma_{11} \: 2 a_{1} a_{2} + \gamma_{20} \: 2 a_{1}^{2} + \gamma_{02} \: 2 a_{2}^{2}.\nonumber 
\end{align}
Having introduced $r$th-order directional derivatives for barycentric polynomials \eqref{eq.12} and \eqref{eq.13}, we may now proceed to the construction of B-spline bases that are $r$th-order continuous throughout.  

\section{Enforcing Continuity}
We now will construct B-spline bases that enforce $r$th-order continuity throughout. We will show how to enforce zeroth- and first-order continuity for second-order polynomial functions defined on the $\left(x,y\right)$ plane. It is left to the reader to generalize to higher variate functions and higher order polynomials.

Say, we have the second-order polynomial:
\begin{equation}
	\label{eq.22} 
		p\!\left(x,y\right) = \alpha_{1} + \alpha_{2} \: x + \alpha_{3} \:  y + \alpha_{4} \:  x y + \alpha_{5} \:  x^{2} + \alpha_{6} \:  y^{2}.
\end{equation}
Then we may make a change of variable from \eqref{eq.22} to a polynomial which takes as its arguments the barycentric coordinates relative to the sides of the triangle $T_{1}$:
\begin{equation}
	\label{eq.23} 
 p_{1}\!\left(b_{1}, b_{2}, b_{3}\right) = \gamma_{00} \:  b_{3}^{2} + \gamma_{10} \:  b_{1}   b_{3} + \gamma_{01} \:  b_{2}   b_{3} + \gamma_{11} \:  b_{1}   b_{2} + \gamma_{20} \:  b_{1}^{2} + \gamma_{02} \:  b_{2}^{2}.
\end{equation}
So, if $\left(x,y\right)\in T_{1}$, then we have $0\leq b_{i} \leq 1$, for $i = 1, 2, 3$. Likewise, may make a change of variable from \eqref{eq.22} to a polynomial which takes as its arguments the barycentric coordinates relative to the sides of the triangle $T_{2}$:
\begin{equation}
	\label{eq.24} 
 p_{2}\!\left(\tilde{b}_{1}, \tilde{b}_{2}, \tilde{b}_{3}\right) = \tilde{\gamma}_{00} \:  \tilde{b}_{3}^{2} + \tilde{\gamma}_{10} \:  \tilde{b}_{1}   \tilde{b}_{3} + \tilde{\gamma}_{01} \:  \tilde{b}_{2}   \tilde{b}_{3} + \tilde{\gamma}_{11} \:  \tilde{b}_{1}   \tilde{b}_{2} + \tilde{\gamma}_{20} \:  \tilde{b}_{1}^{2} + \tilde{\gamma}_{02} \:  \tilde{b}_{2}^{2}.
\end{equation}
So, if $\left(x,y\right)\in T_{2}$, then $0\leq \tilde{b}_{i} \leq 1$, for $i = 1, 2, 3$.

The transformed polynomials \eqref{eq.23} and \eqref{eq.24} are valid on the whole $x,y$-plane. However, we will constrain the polynomials $p_{1}$ and $p_{2}$ to the triangles $T_{1}$ and $T_{2}$, respectively. The basis for the unconnected piecewise polynomials $p_{i}$ of the triangulation $T_{i}$, then may be given as:    
\begin{equation}
	\label{eq.25} 
	B = \left[
\begin{array}{cccccccccccc}
	b_{3}^{2} & b_{1} b_{3} & b_{2} b_{3} & b_{1} b_{2} & b_{1}^{2} & b_{2}^{2} & 0 & 0 & 0 & 0 & 0 & 0 \\
	0 & 0 & 0 & 0 & 0 & 0 & \tilde{b}_{3}^{2} & \tilde{b}_{1} \tilde{b}_{3} & \tilde{b}_{2} \tilde{b}_{3} & \tilde{b}_{1} \tilde{b}_{2} & \tilde{b}_{1}^{2} & \tilde{b}_{2}^{2} 
\end{array}
\right],
\end{equation}
where the columns of the basis $B$ correspond, respectively, with the coefficients
\[
\left(
\begin{array}{cccccccccccc}
	\gamma_{00} & \gamma_{10} & \gamma_{01} & \gamma_{11} & \gamma_{20} & \gamma_{02} & \tilde{\gamma}_{00} & \tilde{\gamma}_{10} & \tilde{\gamma}_{01} & \tilde{\gamma}_{11} & \tilde{\gamma}_{20} & \tilde{\gamma}_{02} 
\end{array}
\right).
\]
If $\left(x,y\right) \in T_{1}$, then the corresponding barycentric coordinates, found by way of \eqref{eq.4} and \eqref{eq.10}, can be fed into the first row. But if $\left(x,y\right) \in T_{2}$, then the corresponding barycentric coordinates are fed into the second row.

Now, if we look at the triangulation in Figure~\ref{fig:triangulation} and equations \eqref{eq.8} and \eqref{eq.9}, then we may see that all points $\left(b_{1}, 0, b_{3}\right)$ on the boundary of $T_{1}$ are equivalent to the points $\left(\tilde{b}_{1}, \tilde{b}_{2}, 0\right)$ on the boundary of $T_{2}$, if $b_{1} = \tilde{b}_{1}$ and $b_{3} = \tilde{b}_{2}$. Stated differently, in order for the polynomials \eqref{eq.23} and \eqref{eq.24} to be connected at their shared boundary, we must have that the zeroth-order directional derivative of these polynomials generate the same values for $q_{1}$ and  $q_{2}$, where
\begin{equation}
	\label{eq.26}
	q_{1} = b_{1} =\tilde{b}_{1}, \qquad  q_{2} = b_{3} =\tilde{b}_{2}, \qquad b_{2} = \tilde{b}_{3} = 0.
\end{equation}
By substituting the values \eqref{eq.26} into \eqref{eq.19}, we find for the piecewise polynomial defined on $T_{1}$:
\begin{equation}
	\label{eq.27}
	 D_{\boldsymbol{a}}^{\left(0\right)} p_{1}\!\left(q_{1}, 0, q_{2}\right) = \gamma_{00} \:  q_{2}^{2} + \gamma_{10} \:  q_{1}   q_{2}  \gamma_{20} \:  q_{1}^{2} + \gamma_{02}.
\end{equation}
And for the piecewise polynomial defined on $T_{2}$, we find:
\begin{equation}
	\label{eq.28}
	 D_{\tilde{\boldsymbol{a}}}^{\left(0\right)} p\!\left(q_{1}, q_{2}, 0\right) = \gamma_{11} \:  q_{1}   q_{2} + \gamma_{20} \:  q_{1}^{2} + \gamma_{02} \:  q_{2}^{2}.
\end{equation}
Substituting \eqref{eq.27} and \eqref{eq.28} into \eqref{eq.25}, we obtain the constraint matrix, say, $Q$:
\begin{equation}
	\label{eq.29} 
	Q = \left[
\begin{array}{cccccccccccc}
	q_{2}^{2} & q_{1} q_{2} & 0 & 0 & q_{1}^{2} & 0 & 0 & 0 & 0 & 0 & 0 & 0 \\
	0 & 0 & 0 & 0 & 0 & 0 & 0 & 0 & 0 & q_{1} q_{2} & q_{1}^{2} & q_{2}^{2} 
\end{array}
\right].
\end{equation}
It follows that zeroth-order continuity may be enforced by merging Columns 1 and 12, Columns 2 and 10, and Columns 5 and 11 of the bases \eqref{eq.25}. In other words, if we constrain the coefficients in \eqref{eq.23} and \eqref{eq.24} to adhere to:
\begin{equation} 
	\label{eq.30}
	\gamma_{20} = \tilde{\gamma}_{20}, \qquad  \gamma_{00} = \tilde{\gamma}_{02}, \qquad  \gamma_{10} = \tilde{\gamma}_{11}.
\end{equation}
Substituting the constraints \eqref{eq.30} into, say, \eqref{eq.24}, and rearranging the terms in both \eqref{eq.23} and \eqref{eq.24} so that equal coefficients are placed beneath each other, we obtain the barycentric polynomials
\begin{equation}
	\label{eq.30b} 
 p_{1}\!\left(b_{1}, b_{2}, b_{3}\right) = \gamma_{00} \:  b_{3}^{2} + \gamma_{10} \:  b_{1}   b_{3} + \gamma_{20} \:  b_{1}^{2} + \gamma_{01} \:  b_{2}   b_{3} + \gamma_{11} \:  b_{1}   b_{2}  + \gamma_{02} \:  b_{2}^{2}.
\end{equation}
and
\begin{equation}
	\label{eq.31} 
 p_{2}\!\left(\tilde{b}_{1}, \tilde{b}_{2}, \tilde{b}_{3}\right) = \gamma_{00} \:  \tilde{b}_{2}^{2} + \gamma_{10} \:  \tilde{b}_{1}   \tilde{b}_{2} + \gamma_{20} \:  \tilde{b}_{1}^{2} + \tilde{\gamma}_{00} \:  \tilde{b}_{3}^{2} + \tilde{\gamma}_{10} \:  \tilde{b}_{1}   \tilde{b}_{3} + \tilde{\gamma}_{01} \:  \tilde{b}_{2}   \tilde{b}_{3} .
\end{equation}

So, by way of \eqref{eq.30b} and \eqref{eq.31}, we find the basis that enforces zeroth-order continuity, that is, $C^{0}$, to be:
\begin{equation}
	\label{eq.32} 
	B = \left[
\begin{array}{ccccccccc}
	b_{3}^{2} & b_{1} b_{3} & b_{1}^{2} & b_{1} b_{2} &  b_{2} b_{3} & b_{2}^{2} & 0 & 0 & 0 \\
	\tilde{b}_{2}^{2} & \tilde{b}_{1} \tilde{b}_{2} & \tilde{b}_{1}^{2}  & 0 & 0 & 0 & \tilde{b}_{3}^{2} & \tilde{b}_{1} \tilde{b}_{3} & \tilde{b}_{2} \tilde{b}_{3} 
\end{array}
\right],
\end{equation}
where the columns of the basis $B$ correspond, respectively, with the coefficients
\[
\left(
\begin{array}{ccccccccc}
	 \gamma_{00}& \gamma_{10}& \gamma_{20}& \gamma_{01}& \gamma_{11}& \gamma_{02}& \tilde{\gamma}_{00}& \tilde{\gamma}_{10}& \tilde{\gamma}_{01} 
\end{array}
\right).
\]

We now proceed to find the basis $B$ that enforces first-order continuity. For the triangles $T_{1}$ and $T_{2}$, any directional vector $\boldsymbol{u}$ which is non-parallel to the boundary shared by $T_{1}$ and $T_{2}$ will suffice:
\begin{equation}
	\label{eq.33} 
	\boldsymbol{u} = \left(1, 0\right)^{T} = \boldsymbol{v}_{2}, \qquad \text{or, equivalently,} \qquad  \boldsymbol{a} = \left(0, 1, 0\right)^{T},
\end{equation}
where we have used \eqref{eq.14} to express $\boldsymbol{u}$ in terms of barycentric coordinates. By substituting both $\boldsymbol{u}$ and the vertices of triangle $T_{2}$ into \eqref{eq.10} and \eqref{eq.4}, we may obtain an equivalent directional vector in the barycentric coordinates of $T_{2}$: 
\begin{equation}
	\label{eq.34} 
	\tilde{\boldsymbol{a}} = \left(1, 1, -1\right)^{T}.
\end{equation}
It may be checked, by way of \eqref{eq.14}, that
\[
	\tilde{\boldsymbol{v}}_{1} + \tilde{\boldsymbol{v}}_{2} - \tilde{\boldsymbol{v}}_{3} = \left(1, 0\right)^{T} = \boldsymbol{u}.  
\]
If we take the directional vector $\boldsymbol{a}$ from \eqref{eq.33} and substitute it in \eqref{eq.20}, we obtain:
\begin{equation}
	\label{eq.35}
	 D_{\boldsymbol{a}}^{\left(1\right)} p_{1}\!\left(b_{1}, b_{2}, b_{3}\right) = \gamma_{01} \: b_{3} + \gamma_{11} \: b_{1} +  \gamma_{02} \: 2 b_{2}.
\end{equation}
The directional derivative of \eqref{eq.31} is:
\begin{align}
	\label{eq.36}
	 D_{\tilde{\boldsymbol{a}}}^{\left(1\right)} p_{2}\!\left(\tilde{b}_{1}, \tilde{b}_{2}, \tilde{b}_{3}\right) &=  \tilde{\boldsymbol{a}}^{T} \left(
\begin{array}{c}
		\gamma_{10} \:  \tilde{b}_{2} + \gamma_{20} \:  2 \tilde{b}_{1} + \tilde{\gamma}_{10} \:  \tilde{b}_{3}\\
		\gamma_{00} \:  2 \tilde{b}_{2} + \gamma_{10} \:  \tilde{b}_{1} + \tilde{\gamma}_{01} \:  \tilde{b}_{3}\\
		\tilde{\gamma}_{00} \:  2 \tilde{b}_{3} + \tilde{\gamma}_{10} \:  \tilde{b}_{1} + \tilde{\gamma}_{01} \:  \tilde{b}_{2}
\end{array}
\right) \nonumber \\
	   \\
	  &= \gamma_{00} \: 2 \tilde{a}_{2} \tilde{b}_{2} +  \gamma_{10} \left(\tilde{a}_{1}\tilde{b}_{2}+\tilde{a}_{2}\tilde{b}_{1}\right) + \gamma_{20} \: 2 \tilde{a}_{1}\tilde{b}_{1} \nonumber \\
	  & \left.\qquad\right. + \tilde{\gamma}_{00} \:  2 \tilde{a}_{3} \tilde{b}_{3} + \tilde{\gamma}_{10} \left(\tilde{a}_{1}\tilde{b}_{3} + \tilde{a}_{3}\tilde{b}_{1}\right) + \tilde{\gamma}_{01}  \left(\tilde{a}_{2}\tilde{b}_{3} + \tilde{a}_{3}\tilde{b}_{2}\right).\nonumber 
\end{align}
If we substitute the directional vector $\tilde{\boldsymbol{a}}$ from \eqref{eq.34} in \eqref{eq.36}, we obtain:
\begin{align}
	\label{eq.37}
	 D_{\tilde{\boldsymbol{a}}}^{\left(1\right)} p_{2}\!\left(\tilde{b}_{1}, \tilde{b}_{2}, \tilde{b}_{3}\right) &=\gamma_{00} \: 2  \tilde{b}_{2} +  \gamma_{10} \left(\tilde{b}_{2}+\tilde{b}_{1}\right) + \gamma_{20} \: 2 \tilde{b}_{1} \nonumber \\
	  & \left.\qquad\right. - \tilde{\gamma}_{00} \:  2 \tilde{b}_{3}  + \tilde{\gamma}_{10} \left(\tilde{b}_{3} - \tilde{b}_{1}\right) + \tilde{\gamma}_{01} \left(\tilde{b}_{3} - \tilde{b}_{2}\right).
\end{align}
Further, substituting \eqref{eq.26} into \eqref{eq.35} and \eqref{eq.37}, we obtain, respectively,
\begin{equation}
	\label{eq.38}
	 D_{\boldsymbol{a}}^{\left(1\right)} p_{1}\!\left(q_{1}, 0, q_{2}\right) = \gamma_{01} \: q_{2} + \gamma_{11} \: q_{1} 
\end{equation}
and
\begin{equation}
	\label{eq.39}
	 D_{\tilde{\boldsymbol{a}}}^{\left(1\right)} p_{2}\!\left(q_{1}, q_{2}, 0\right) = \gamma_{00} \: 2 q_{2} + \gamma_{10} \left(q_{1} + q_{2}\right) + \gamma_{20} \: 2 q_{1} - \tilde{\gamma}_{10} q_{1} - \tilde{\gamma}_{01} q_{2}. 
\end{equation}
Equations \eqref{eq.38} and \eqref{eq.39} correspond with the constraint matrix:       
\begin{equation}
	\label{eq.40} 
	Q = \left[
\begin{array}{ccccccccc}
	0 & 0 & 0 & q_{1} & q_{2} & 0 & 0 & 0 & 0 \\
	2 q_{2} & q_{1} + q_{1}  & 2 q_{1}  & 0 & 0 & 0 & 0 & - q_{1} & - q_{2} 
\end{array}
\right].
\end{equation}

Now, any permutation of the columns of $B$ will leave intact any $C^{r}$ constraints already in place. For example, the columns of \eqref{eq.32} enforce $C^{0}$, and any permutation of these columns will still enforce this continuity constraint. So, for some polynomial function defined on a $n$-dimensional hyperplane, we are free to permutate the columns of any constraint matrix $C$, in order to obtain rows in which we are left with terms of the type
\begin{equation}
	\label{eq.41} 
	q_{1}^{k_{1}} q_{2}^{k_{2}} \cdots q_{n}^{k_{n}},
\end{equation}
where $\sum_{i=1}^{n} k_{i} = d - r$, and where $d$ is the order of the polynomial function and $r$ is the order of the directional derivative or, equivalently, the order of the continuity constraint we wish to enforce. 

If we thus reduce the rows of our constraint matrix $Q$, then we can instantly see how to enforce our continuity constraint. Just merge the permutated columns having identical elements. For example, in \eqref{eq.29} we have the situation that no permutation matrix is needed, seeing that the elements of the constraint matrix $Q$ are already in the form \eqref{eq.41}. This holds generally for zeroth-order constraint matrices $Q$. 

In what follows, we give the steps needed to derive a permutation matrix that will deliver us \eqref{eq.40} in the desired form. First, we observe in \eqref{eq.40} that the columns 6 and 7 of basis \eqref{eq.32} play no role in the enforcement of the first-order continuity. So, for the moment we will set these columns aside. This leaves us with seven non-zero columns. Now, the dimensionality of our original polynomial \eqref{eq.22} is $n = 2$, the order of this polynomial is $d = 2$, and the continuity constraint we wish to enforce is $r = 1$.  So, we wish to find that permutation matrix that only leaves us with the terms, \eqref{eq.41}: 
\[
	\left\{q_{1}, q_{2}\right\},
\]
in every row of \eqref{eq.40}. 

If we transpose \eqref{eq.40}, then we may designate for each row and each distinct element in the set $\left\{q_{1}, q_{2}\right\}$ a separate column:
\begin{equation}
	\label{eq.42} 
	\left[
\begin{array}{cc}
	0 & 2 q_{2}  \\
	0 & q_{1} + q_{2} \\
	0 & 2 q_{1} \\
	q_{1} & 0 \\
	q_{2} & 0 \\
	0 & - q_{1}\\
	0 & - q_{2}
\end{array}
\right] \longrightarrow  \left[
\begin{array}{cccc}
	0 & 0 & 0 & 2 \\
	0 & 0 & 1 & 1 \\
	0 & 0 & 2 & 0 \\
	1 & 0 & 0 & 0 \\
	0 & 1 & 0 & 0 \\
	0 & 0 & -1 & 0 \\
	0 & 0 & 0 & -1
\end{array}
\right] 
\end{equation}
Then we append an $7 \times 7$ identity matrix to the right-hand side of \eqref{eq.42}:   
\begin{equation}
	\label{eq.43} 
	\left[\begin{array}{ccccccccccc}
	0 & 0 & 0 & 2 & 1 & 0 & 0 & 0 & 0 & 0 & 0\\
	0 & 0 & 1 & 1 & 0 & 1 & 0 & 0 & 0 & 0 & 0\\
	0 & 0 & 2 & 0 & 0 & 0 & 1 & 0 & 0 & 0 & 0\\
	1 & 0 & 0 & 0 & 0 & 0 & 0 & 1 & 0 & 0 & 0\\
	0 & 1 & 0 & 0 & 0 & 0 & 0 & 0 & 1 & 0 & 0\\
	0 & 0 & -1 & 0 & 0 & 0 & 0 & 0 & 0 & 1 & 0\\
	0 & 0 & 0 & -1 & 0 & 0 & 0 & 0 & 0 & 0 & 1
\end{array}
\right] 
\end{equation}
If we row reduce \eqref{eq.43}, by way of Gaussian elimination \cite{Lay00}, then we get:
\begin{equation}
	\label{eq.44} 
	\left[\begin{array}{ccccccccccc}
	1 & 0 & 0 & 0 & 0 & 0 & 0 & 1 & 0 & 0 & 0\\
	0 & 1 & 0 & 0 & 0 & 0 & 0 & 0 & 1 & 0 & 0\\
	0 & 0 & 1 & 0 & 0 & 0 & 0 & 0 & 0 &-1 & 0\\
	0 & 0 & 0 & 1 & 0 & 0 & 0 & 0 & 0 & 0 &-1\\
	0 & 0 & 0 & 0 & 1 & 0 & 0 & 0 & 0 & 0 & 2\\
	0 & 0 & 0 & 0 & 0 & 1 & 0 & 0 & 0 & 1 & 1\\
	0 & 0 & 0 & 0 & 0 & 0 & 1 & 0 & 0 & 2 & 0
\end{array}
\right] 
\end{equation}
By dropping the first four columns of \eqref{eq.44}, we are left with the transposed permutation matrix:
\begin{equation}
	\label{eq.46} 
	\boldsymbol{P}^{T} = \left[\begin{array}{ccccccc}
	 0 & 0 & 0 & 1 & 0 & 0 & 0\\
	 0 & 0 & 0 & 0 & 1 & 0 & 0\\
	 0 & 0 & 0 & 0 & 0 &-1 & 0\\
   0 & 0 & 0 & 0 & 0 & 0 &-1\\
	 1 & 0 & 0 & 0 & 0 & 0 & 2\\
	 0 & 1 & 0 & 0 & 0 & 1 & 1\\
	 0 & 0 & 1 & 0 & 0 & 2 & 0
\end{array}
\right] 
\end{equation}   
The non-zero columns of \eqref{eq.40} are
\begin{equation}
	\label{eq.47} 
	Q_{1} = \left[
\begin{array}{ccccccc}
	0 & 0 & 0 & q_{1} & q_{2} & 0 & 0 \\
	2 q_{2} & q_{1} + q_{1}  & 2 q_{1}  & 0 & 0 & - q_{1} & - q_{2} 
\end{array}
\right].
\end{equation}   
It may be checked that multiplying $Q_{1}$ with the permutation matrix $\boldsymbol{P}$, brings the non-zero columns of \eqref{eq.40} in the desired form:
\begin{equation}
	\label{eq.48} 
	Q_{1} \boldsymbol{P}  = \left[
\begin{array}{ccccccc}
	q_{1} & q_{2} & 0      & 0 		 & 0 & 0 & 0 \\
	0 	  & 0  	  & q_{1}  & q_{2} & 0 & 0 & 0 
\end{array}
\right].
\end{equation}   
Now that we have found the needed permutation matrix $\boldsymbol{P}$, we may proceed to the construction of a first-order continuous B-spline. 

We select the columns in \eqref{eq.32} which are non-zero \eqref{eq.40}:
\begin{equation}
	\label{eq.49} 
	B_{1} = \left[
\begin{array}{ccccccc}
	b_{3}^{2} & b_{1} b_{3} & b_{1}^{2} & b_{1} b_{2} &  b_{2} b_{3} &  0 & 0 \\
	\tilde{b}_{2}^{2} & \tilde{b}_{1} \tilde{b}_{2} & \tilde{b}_{1}^{2}  & 0 & 0 & \tilde{b}_{1} \tilde{b}_{3} & \tilde{b}_{2} \tilde{b}_{3} 
\end{array}
\right]
\end{equation}
The selected columns of \eqref{eq.32} are then multiplied with the permutation matrix $\boldsymbol{P}$:
\begin{align}
	\label{eq.50} 
	B_{2} &= B_{1} \boldsymbol{P} \nonumber\\
	\\
	&= \left[
\begin{array}{ccccccc}
	b_{1} b_{2} & b_{2} b_{3} &	0	&	0	&	b_{3}^{2} & b_{1} b_{3} & b_{1}^{2} \\
	0	&	0	& - \tilde{b}_{1} \tilde{b}_{3} & - \tilde{b}_{2} \tilde{b}_{3}	& \tilde{b}_{2}^{2} + 2 \tilde{b}_{2} \tilde{b}_{3}  & \tilde{b}_{1} \tilde{b}_{2} + \tilde{b}_{1} \tilde{b}_{3} + \tilde{b}_{2} \tilde{b}_{3}  & \tilde{b}_{1}^{2} + 2 \tilde{b}_{1} \tilde{b}_{3}
\end{array}
\right]. \nonumber
\end{align}
We then, because of \eqref{eq.48}, merge Columns 1 and 3, and Columns 2 and 4:
\begin{equation}
	\label{eq.51} 
	B_{3} = \left[
\begin{array}{ccccc}
	b_{1} b_{2} & b_{2} b_{3} &	b_{3}^{2} & b_{1} b_{3} & b_{1}^{2} \\
	- \tilde{b}_{1} \tilde{b}_{3} & - \tilde{b}_{2} \tilde{b}_{3}	& \tilde{b}_{2}^{2} + 2 \tilde{b}_{2} \tilde{b}_{3}  & \tilde{b}_{1} \tilde{b}_{2} + \tilde{b}_{1} \tilde{b}_{3} + \tilde{b}_{2} \tilde{b}_{3}  & \tilde{b}_{1}^{2} + 2 \tilde{b}_{1} \tilde{b}_{3}
	\end{array}
\right]. 
\end{equation}
Finally, we add the columns in \eqref{eq.32} which were dropped in \eqref{eq.49}: 
\begin{equation}
	\label{eq.52} 
	B_{4} = \left[
\begin{array}{ccccccc}
	b_{1} b_{2} & b_{2} b_{3} &	b_{3}^{2} & b_{1} b_{3} & b_{1}^{2} & b_{2}^{2} & 0\\
	- \tilde{b}_{1} \tilde{b}_{3} & - \tilde{b}_{2} \tilde{b}_{3}	& \tilde{b}_{2}^{2} + 2 \tilde{b}_{2} \tilde{b}_{3}  & \tilde{b}_{1} \tilde{b}_{2} + \tilde{b}_{1} \tilde{b}_{3} + \tilde{b}_{2} \tilde{b}_{3}  & \tilde{b}_{1}^{2} + 2 \tilde{b}_{1} \tilde{b}_{3} & 0 & \tilde{b}_{3}^{2} 
	\end{array}
\right]. 
\end{equation}  
Seeing that the columns of this basis are a linear combination of the columns of the basis \eqref{eq.32}, which was constructed to be zeroth-order continuous throughout, that is, $C^{0}$, we have that the basis \eqref{eq.52} is first-order continuous throughout, that is, $C^{1}$.

In order to construct a B-spline basis that is second-order continuous throughout, that is, $C^{2}$, we simply repeat the steps that led us to \eqref{eq.52}. 

First, we take the second-order directional derivative of \eqref{eq.52}. By doing so, we may identify and put aside the columns in \eqref{eq.52} which do not participate in the second-order continuity constraint, those being the zero-columns of the second-order directional derivative. This then gives us a reduced first-order explicit base. We then proceed to construct a permutation matrix $\boldsymbol{P}$ for the remaining, that is, non-zero, columns of the second-order directional derivative. If we multiply the reduced first-order explicit base with this permutation matrix $\boldsymbol{P}$, and adding again those columns of \eqref{eq.52} that were put aside, then we obtain the second order B-spline basis that is second-order continuous throughout, that is, $C^{2}$. 

By an iteration of these steps, we may arrive at the explicit B-spline basis that is $C^{r}$, for $r \geq 3$. If we expand the triangulation in Figure~\ref{fig:triangulation}, in that we add more triangles, we must, for a given $C^{r}$, repeat these steps for every pair of triangles which share a side. 

Note that for $n$-variate B-splines, which have $n$-variate triangulations, all tetrahedra which share $\left(n-1\right)$-dimensional sides are considered to be connected.

\section{A Simple B-Spline Analysis}
Say, we have a triangulation as in Figure~\ref{fig:triangulation}, a polynomial of degree $d = 2$ and $r = 1$ continuity throughout, then the corresponding B-spline base is \eqref{eq.52}:
\begin{equation}
	\label{eq.53} 
	B = \left[
\begin{array}{ccccccc}
	b_{1} b_{2} & b_{2} b_{3} &	b_{3}^{2} & b_{1} b_{3} & b_{1}^{2} & b_{2}^{2} & 0\\
	- \tilde{b}_{1} \tilde{b}_{3} & - \tilde{b}_{2} \tilde{b}_{3}	& \tilde{b}_{2}^{2} + 2 \tilde{b}_{2} \tilde{b}_{3}  & \tilde{b}_{1} \tilde{b}_{2} + \tilde{b}_{1} \tilde{b}_{3} + \tilde{b}_{2} \tilde{b}_{3}  & \tilde{b}_{1}^{2} + 2 \tilde{b}_{1} \tilde{b}_{3} & 0 & \tilde{b}_{3}^{2} 
	\end{array}
\right].
\end{equation}  
The first and second row of $B$ correspond, respectively, with $T_{1}$ and $T_{2}$ of Figure~\ref{fig:triangulation}. 

If we have a small data set of $n = 5$ observations $\left(x_{i}, y_{i}, z_{i}\right)$: 
\begin{align}
	\label{eq.54} 
	\left(x_{1}, y_{1}, z_{1}\right) &= \left(0.2, 0.1, 1.0\right), \nonumber \\
	\left(x_{2}, y_{2}, z_{2}\right) &= \left(0.2, 0.7, 3.0\right), \nonumber \\
	\left(x_{3}, y_{3}, z_{3}\right) &= \left(0.1, 0.3, 2.0\right),		 \\
	\left(x_{4}, y_{4}, z_{4}\right) &= \left(0.5, 0.1, 1.0\right),	\nonumber \\
	\left(x_{5}, y_{5}, z_{5}\right) &= \left(0.7, 0.8, 4.0\right). \nonumber
\end{align}  
Then we have that the input values, both in Cartesian and barycentric coordinates, \eqref{eq.10} and \eqref{eq.4}, are given as:
\begin{align}
	\label{eq.55} 
	\left(x_{1}, y_{1}\right) &= \left(0.2, 0.1\right) \in T_{1}, \qquad &\left(b_{11}, b_{12}, b_{13}\right) = \left(0.8, 0.1, 0.1\right), \nonumber \\
	\left(x_{2}, y_{2}\right) &= \left(0.2, 0.7\right) \in T_{2}, \qquad &\left(\tilde{b}_{21}, \tilde{b}_{22}, \tilde{b}_{23}\right) = \left(0.3, 0.2, 0.5\right), \nonumber \\
	\left(x_{3}, y_{3}\right) &= \left(0.1, 0.3\right) \in T_{2}, \qquad &\left(\tilde{b}_{31}, \tilde{b}_{32}, \tilde{b}_{33}\right) = \left(0.7, 0.1, 0.2\right),  \\
	\left(x_{4}, y_{4}\right) &= \left(0.5, 0.1\right) \in T_{1}, \qquad &\left(b_{41}, b_{42}, b_{43}\right) = \left(0.5, 0.4, 0.1\right), \nonumber \\
	\left(x_{5}, y_{5}\right) &= \left(0.7, 0.8\right) \in T_{2}, \qquad &\left(\tilde{b}_{51}, \tilde{b}_{52}, \tilde{b}_{53}\right) = \left(0.2, 0.7, 0.1\right), \nonumber
\end{align}  
The vector with corresponding output values is
\begin{equation}
	\label{eq.56} 
	\boldsymbol{z} = \left(
\begin{array}{ccccc}
	1.0 & 3.0 & 2.0 & 1.0 & 4.0
  \end{array}
\right)^{T}. 
\end{equation} 
The points $\left(b_{11}, b_{21}, b_{31}\right)$ and  $\left(b_{14}, b_{24}, b_{34}\right)$ are assigned to the first partitioning, or, equivalently, to the first row of $B$. Likewise, $\left(b_{12}, b_{22}, b_{32}\right)$, $\left(b_{13}, b_{23}, b_{33}\right)$, and $\left(b_{15}, b_{25}, b_{35}\right)$ are assigned to the second partitioning, or, equivalently, the second row of $B$:
\begin{align}
	\label{eq.57} 
	B &= \left[
\begin{array}{ccccccc}
	b_{11} b_{12} & b_{12} b_{13} &	b_{13}^{2} & b_{11} b_{13} & b_{11}^{2} & b_{12}^{2} & 0\\
	- \tilde{b}_{21} \tilde{b}_{23} & - \tilde{b}_{22} \tilde{b}_{23}	& \tilde{b}_{22}^{2} + 2 \tilde{b}_{22} \tilde{b}_{23}  & \tilde{b}_{21} \tilde{b}_{22} + \tilde{b}_{21} \tilde{b}_{23} + \tilde{b}_{22} \tilde{b}_{23}  & \tilde{b}_{21}^{2} + 2 \tilde{b}_{21} \tilde{b}_{23} & 0 & \tilde{b}_{23}^{2} \\
	- \tilde{b}_{31} \tilde{b}_{33} & - \tilde{b}_{32} \tilde{b}_{33}	& \tilde{b}_{32}^{2} + 2 \tilde{b}_{32} \tilde{b}_{33}  & \tilde{b}_{31} \tilde{b}_{32} + \tilde{b}_{31} \tilde{b}_{33} + \tilde{b}_{32} \tilde{b}_{33}  & \tilde{b}_{31}^{2} + 2 \tilde{b}_{31} \tilde{b}_{33} & 0 & \tilde{b}_{33}^{2} \\
	b_{41} b_{42} & b_{42} b_{43} &	b_{43}^{2} & b_{41} b_{43} & b_{41}^{2} & b_{42}^{2} & 0\\
	- \tilde{b}_{51} \tilde{b}_{53} & - \tilde{b}_{52} \tilde{b}_{53}	& \tilde{b}_{52}^{2} + 2 \tilde{b}_{52} \tilde{b}_{53}  & \tilde{b}_{51} \tilde{b}_{52} + \tilde{b}_{51} \tilde{b}_{53} + \tilde{b}_{52} \tilde{b}_{53}  & \tilde{b}_{51}^{2} + 2 \tilde{b}_{51} \tilde{b}_{53} & 0 & \tilde{b}_{53}^{2} 
	\end{array}
\right] \nonumber \\
\\
&= \left[
\begin{array}{ccccccc}
	0.08  &  0.01 &	0.01 & 0.08 & 0.64 & 0.01 & 0 \\
	-0.15 & -0.10 &	0.24 & 0.31 & 0.39 & 0 		& 0.25 \\
	-0.14 & -0.02 &	0.05 & 0.23 & 0.77 & 0 		& 0.04 \\
	0.20  &  0.04 &	0.01 & 0.05 & 0.25 & 0.16 & 0 \\
	-0.02 & -0.07 &	0.63 & 0.23 & 0.08 & 0 		& 0.01 \\
	\end{array}
\right]. \nonumber
\end{align}  
The unknown regression coefficients $\boldsymbol{\gamma}$ of the B-spline may be found as the solution of a simple regression problem, \eqref{eq.56} and \eqref{eq.57}:
\begin{equation}
	\label{eq.58} 
	\boldsymbol{\gamma} = \left(B^{T}B\right)^{-1}B^{T}\boldsymbol{z}.
\end{equation} 

The B-spline interpolation estimate, for some new set of barycentric coordinates, 
\[
	\left(b_{1}, b_{2}, b_{3}\right),
\]
then may be found by substituting these coordinates in the appropriate row of \eqref{eq.53}. This will result in some $7 \times 1$ vector, say, $\boldsymbol{b}$. The interpolation estimate $\hat{z}$ then may be found as the inner product, \eqref{eq.58}: 
\begin{equation}
	\label{eq.59} 
	\hat{z} = \boldsymbol{b}^{T}\boldsymbol{\gamma}.
\end{equation} 
This concludes our short outline on how to set up a simple B-spline analysis.

\section{Discussion}
We have introduced here a direct method to construct explicit B-spline bases. In this direct method, for a given $C^{r}$, the computational burden lies with the row reduction, needed for the construction of the permutation matrices $\boldsymbol{P}$, for each pair of connecting tetrahedra. How this direct method compares with the indirect Lagrangian method of Awanou, \cite{Awanou03}, will be the subject of a future paper. 

Finally, we refer the interested reader to \cite{vanErp14}, where a method for the construction of explicit C-spline bases is given. C-splines are piecewise polynomials which are $C^{r}$ continuous throughout, defined on adjacent Cartesian, as opposed to barycentric, coordinate systems. C-splines constitute a generalization of the ordinary linear regression model of statistics; in that the latter are a special case of the former.


\begin{thebibliography}{9}                                                                                                %


\bibitem {Awanou03} 
Awanou G.M. (2003), \textit{Energy methods in 3D Spline approximations of the Navier-Stokes equations}, thesis work under the direction of Ming-Jun Lai; Athens, Georgia.

\bibitem {Lay00}
Lay, D.C. (2000), \textit{Linear algebra and its applications}, Addison-Wesley Publishing Company; 2nd edn update.
\end{thebibliography}
\end{document}